\documentclass[a4paper,14pt]{article}
\usepackage{amsmath,amsfonts,amssymb,amscd,amsthm,epsf}
\usepackage{times}
\makeatletter\@ifundefined{email}{\def\email#1{}}{}\makeatother%
\usepackage[T1]{fontenc}
\usepackage{url}
\usepackage {color}
  \definecolor{Light}{gray}{.80}
  \definecolor{Dark}{gray}{.20}
  \usepackage{boxedminipage}
  \usepackage{shadow}
  \usepackage{colordvi}
  \usepackage{graphics}

\def\N{\mathbf{N}} 
\def\Z{\mathbf{Z}} 

\title{GEOMETRIE ET COGNITION:\\L'EXEMPLE DU CONTINU\\  \small CERISY SEPTEMBRE 2006}\large
\author{B. Teissier\\
Equipe `` G\'eom\'etrie et Dynamique'',\\
Institut Math\'ematique de Jussieu,\\ UMR 7586 du CNRS\\
175 Rue du Chevaleret\\
F-75013 Paris, France.\\ email: teissier@math.jussieu.fr }
\date{\today}
\begin{document}

\maketitle
{\small \centerline{\textbf{Summary}}
In this paper I propose the idea to establish a clear distinction between the foundations of truth and the foundations of meaning in Mathematics. I explore on the most basic example, the mathematical line, the possibility that the foundations of its meaning are provided by a protomathematical object resulting from the identification by our perceptual system of the visual line and the vestibular line, a point of view suggested by recent results of neurophysiology.}

\medskip

\section {INTRODUCTION}

\par\medskip\noindent
{\it Quand nous pensons une pens\'ee, la signification de cette pens\'ee est la forme du processus neuro-physiologique sous-jacent.}\par\hfill Bernhard Riemann, {\it opuscule philosophique}.
\par\medskip\noindent
{\it En r\'esum\'e, pour chaque attitude de mon corps, mon premier doigt d\'etermine un point et c'est cela, et cela seulement qui d\'efinit un point de l'espace.}\par\hfill
Henri Poincar\'e, {\it La Science et l'hypoth\`ese}, Flammarion.
\par\bigskip
\par\medskip\noindent
{\it ...et cela para\^\i tra bien beau \`a l'homme instruit, mais d\'eplaira beaucoup aux ignorants}.\par\hfill
Ma\"\i monide, {\it Le guide des \'egar\'es}.
\par\bigskip

Thierry Paul a commenc\'e son expos\'e en nous disant que tout allait bien pour la m\'ecanique quantique: l'accord avec l'exp\'erience est excellent et il n'y a plus de paradoxes g\^enants.\par\noindent
De m\^eme, les math\'ematiques ne se sont jamais si bien port\'ees, les conjectures et probl\`emes tombent les uns apr\`es les autres, les id\'ees foisonnent, de nouveaux domaines naissent.\par
A la diff\'erence de la Physique (qui s'estime probablement fond\'ee {\it de facto} par l'accord avec l'exp\'erience, comme nous l'a rappel\'e Thierry Paul), et ne revient sur ses fondements qu'en cas de crise majeure comme celle des d\'ebuts de la m\'ecanique quantique,  les math\'ematiques, si prosp\`eres qu'elles soient, s'int\'eressent constamment  \`a leurs fondements.\par Comme nous savons, depuis presque un  si\`ecle  l'aspect "philosophie naturelle" cher \`a Riemann, Hermann Weyl, Poincar\'e, Thom, a trop souvent \'et\'e d\'elaiss\'e dans la r\'eflexion sur les fondements au profit de l'aspect axiomatique. Celui-ci fonde les math\'ematiques d'une mani\`ere en apparence inattaquable et leur donne une puissance cr\'eative inimaginable auparavant tout en rendant possible l'oubli de leur enracinement historique dans le r\'eel et en bouleversant la construction de leur signification. On peut dire sans trop exag\'erer, en particulier en ce qui concerne l'enseignement, que les racines "naturelles" des math\'ematiques ont \'et\'e remplac\'ees, apr\`es la crise des g\'eom\'etries non euclidiennes, par une chape de b\'eton devant permettre de construire solidement. Les math\'ematiques n'\'etant plus enracin\'ees dans le r\'eel par nature, le probl\`eme de leurs fondements est apparu avec urgence, mais il est apparu beaucoup plus clairement en ce qui concerne la v\'erit\'e, susceptible d'une rigoureuse analyse syntaxique, qu'en ce qui concerne la signification qui pour beaucoup de scientifiques \'etait trop difficile \`a appr\'ehender dans les brouillards de la "psychologie".\par Une bonne partie de la philosophie des math\'ematiques s'est donc focalis\'ee sur le seul aspect logique (au sens classique, pas celui de LIGC) de cette magnifique entreprise de refondation et sur les probl\`emes math\'ematiques que posent la th\'eorie des ensembles et la th\'eorie de la d\'emonstration.\par Comme l'a rappel\'e aussi Thierry Paul, cette p\'eriode d'interrogation sur les fondements au d\'ebut du XX$^{\hbox{\rm\`eme}}$ si\`ecle a aussi \'et\'e une p\'eriode extr\^emement f\'econde en Math\'ematiques, et l'activit\'e fondationnelle ne repr\'esentait qu'une toute petite partie de l'activit\'e totale. \par Chez les math\'ematiciens l'int\'er\^et pour les fondements est \'evidemment plus fort dans les p\'eriodes o\`u la discipline a besoin d'examiner des objets si nouveaux que son support intuitif ne peut pas suivre (pensez aux fonctions continues non d\'erivables, aux ordinaux) et qu'ils semblent donner mati\`ere \`a  controverse ou \`a paradoxe. Une entreprise de construction formelle et de r\'evision des significations (deux op\'erations dont j'essaie d'\'etudier les relations dans ce texte) est alors lanc\'ee et, en tous cas jusqu'\`a maintenant, les choses rentrent dans l'ordre. Par exemple la signification de l'ensemble de Cantor est d\'esormais per\c cue par un bon \'el\`eve de Licence. Quoi qu'il en soit, il semble que cet int\'er\^et pour les fondements soit assez faible en ce moment. Voir [G3] pour une r\'eflexion stimulante sur les fondements et la logique correspondante.\par
La th\'eorie des ensembles de Zermelo-Frenkel semble suffisante pour donner un statut \`a tous les objets dont nous avons besoin, et la ''logique na\"\i ve'' du math\'ematicien moderne semble suffisante pour \'eviter les erreurs de raisonnement; ce math\'ematicien n'en demande pas plus sauf si le probl\`eme des fondements le pr\'eoccupe.\par Il a peut-\^etre tort, parce que l'enseignement des math\'ematiques semble de plus en plus \'eloign\'e du ''r\'eel'' et priv\'e de signification, et attire de moins en moins d'\'etudiants. \par
La p\'eriode actuelle semble propice pour se pencher sereinement sur deux probl\`emes qui devraient int\'eresser les philosophes et certains math\'ematiciens et auxquels le ''Mod\`ele standard'' construit sur l'axiomatique de Zermelo-Frenkel n'apporte aucune r\'eponse parce qu'ils rel\`event non pas des fondements de la v\'erit\'e mais des fondements de la signification.

\par\medskip\noindent
$\bullet$  Comment se fait-il que le sentiment de comprendre une d\'emonstration semble si \'eloign\'e de la structure logique de la preuve?\par\medskip\noindent
$\bullet$ Comment d\'etermine-t-on parmi la myriade d'\'enonc\'es, ceux qui sont assez int\'eressants pour que l'on essaie de les d\'emontrer?

C'est le moment de citer un des beaux aphorismes de Ren\'e Thom:\par\medskip\noindent
\centerline {\it La limite de la v\'erit\'e n'est pas l'erreur, c'est l'insignifiance.}\par\medskip\noindent
 Si, encourag\'es par cet aphorisme, nous acceptons l'id\'ee qu'il faut fonder non seulement la v\'erit\'e des th\'eor\`emes, mais aussi leur signification, il semble important de d\'evelopper une analyse scientifique des aspects implicites du raisonnement qui font appel par exemple \`a des jugements sur le continu spatial, sa connexit\'e, ses sym\'etries etc., ou \`a des jugements sur l'ensemble des entiers, ses op\'erations et le bon ordre dont il est \'equip\'e dans notre espace mental.\par
 Des progr\`es dans cette direction permettraient peut-\^etre, par un ''retour au sens'', de rendre les math\'ematiques plus accessibles, et il n'est pas impossible qu'ils aident m\^eme des chercheurs confirm\'es par une augmentation de lucidit\'e.\par La phrase de Riemann reproduite au d\'ebut indique une direction dans laquelle il \'etait jusqu'ici difficile d'avancer faute de pr\'ecisions sur les processus neuro-physiologiques humains, et aussi faute d'int\'er\^et (pour ne pas parler de m\'efiance) \`a l'\'egard des sources ''psychologiques'' de la signification de la part des scientifiques, hormis quelques\break exceptions c\'el\`ebres comme Poincar\'e.\par Au d\'ebut les progr\`es des neurosciences ont surtout servi de base scientifique au raidissement id\'eologique \`a la Changeux ({\it cf}. [C] et [C-C]) du r\'eductionnisme discr\`etement envahissant de la fin du XX$^{\hbox{\rm\`eme}}$ si\`ecle. Plus r\'ecemment ces progr\`es ont aussi servi d'argument pour le d\'eveloppement d'une vision plus riche ({\it cf}. [B2], [D]) de la pens\'ee consciente, faisant intervenir entre autres les \'emotions. Cependant si des ''moteurs'' de type \'emotionnel (involontaires, encore une fois, mais plut\^ot conscients) comme le d\'esir jouent sans doute un r\^ole dans la pens\'ee scientifique, ce r\^ole me semble bien difficile \`a cerner, peut-\^etre pas primordial, et je pr\'ef\`ere me concentrer ici sur les pulsions de nature plus ''structurelle''.\par\medskip\noindent
 S'il est indiscutable qu'un fondement de la v\'erit\'e sans d\'eduction rationnelle est de nature religieuse et n'est pas une option pour un scientifique contemporain, il me semble important de distinguer entre ``Fondements (ou sources) de la v\'erit\'e des Math\'ematiques'' et ``Fondements (ou sources) de la signification des Math\'ematiques''.\par   Il est peut-\^etre temps en effet, plus de deux si\`ecles apr\`es Kant, d'admettre une bonne fois que la signification n'est pas r\'eductible \`a la pens\'ee rationnelle et consciente et surtout d'en {\it tirer les cons\'equences} en essayant de d\'ecrire comment et \`a partir de quoi elle se construit, \`a quels besoins elle r\'epond, et ce qu'elle apporte. Ce que la signification apporte n'est pas une preuve de v\'erit\'e, mais me semble tout aussi indispensable \`a la pens\'ee scientifique, et un d\'ebut d'analyse rigoureuse de sa construction est d\'esormais possible.\par Cette analyse s'appuie d'abord sur les progr\`es dans la compr\'ehension biologique de notre perception du monde, de la m\'emoire et des \'ev\`enements inconscients qui\break foisonnent en nous. Cette compr\'ehension est source de nouveaux concepts, et permet d'imaginer de nouveaux modes d'explication tenant compte de notre rapport inconscient au monde et se substituant aux tentatives vaines d'analyser ''rationnellement'' la signification \`a l'aide de mots comme "m\'etaphore" dont le sens demande lui-m\^eme \`a \^etre expliqu\'e au moyen d'autres mots, et ainsi de suite dans une fuite sans fin analogue \`a celle du ``meta'' d\'enonc\'ee par Girard dans [G1]. Voir aussi [T1].\par J'esp\`ere en particulier convaincre le lecteur que la doctrine qui pr\'evaut implicitement dans l'enseignement selon laquelle la signification des objets math\'ematiques se construit toujours parall\`element aux fondements de leur v\'erit\'e est erron\'ee pour des raisons fondamentales, et pas seulement parce qu'elle est contredite par l'histoire des math\'ematiques et l'exp\'erience des chercheurs.\footnote{C\'eder \`a la tentation (de bas niveau au sens d\'efini un peu plus bas) d'identifier les fondements de la signification avec les fondements de la v\'erit\'e conduit aussi, \`a mon avis, \`a la prolif\'eration de faux probl\`emes en philosophie des math\'ematiques.}
 \par\medskip\noindent
 Cela ne signifie pas que les {\it mod\`eles} de la pens\'ee rationnelle ne puissent pas avoir de signification, puisque apr\`es tout ils refl\`etent aussi nos automatismes de pens\'ee et leur construction utilise notre perception de l'espace, parfois d'une mani\`ere assez explicite comme dans la diagrammatique ch\`ere \`a Gilles Ch\^atelet et Charles Alunni, et parfois d'une mani\`ere tr\`es souterraine. Mais pour le moment il me semble que ce n'est justement que d'une signification de {\it leur construction} et  de leurs interactions (qui sont contraintes par leur construction) qu'il s'agit. Cela diff\`ere de la signification propre que peuvent avoir des objets math\'ematiques. La distinction est cependant pr\'ecaire car la nature de certains objets math\'ematiques \'elabor\'es r\'eside pour une bonne part dans leur construction. C'est le cas par exemple des alg\`ebres d'op\'erateurs qu'utilise Girard.\par
Dans son expos\'e, nous avons appris une mani\`ere de contr\^oler la g\'en\'erativit\'e folle des symboles en interpr\'etant dans une alg\`ebre ''mod\'eratrice'' les op\'erations qu'ils sont cens\'es repr\'esenter, d'une mani\`ere qui pr\'eserve (et m\^eme pr\'ecise) leur signification logique. Nous avons aussi appris que quelque chose comme un ``sens'' ou un ''point de vue'' y influait sur la ''v\'erit\'e''.\par
On peut r\^ever que l'\'etude du ''sens cognitif '' des \'enonc\'es, qui est aussi une mani\`ere de mod\'erer la g\'en\'erativit\'e des constructions formelles, rencontrera un jour la G\'eom\'etrie des Interactions. On peut voir une partie du livre [Pe1] comme la description du d\'ebut d'un chemin allant dans cette direction. Pour une riche description d'autres tentatives dans cette direction, je renvoie \`a [B-L] et [C], ainsi qu'aux textes de J.-Y. Girard et G. Longo dans ce volume. Dans une veine diff\'erente, je recommande vivement le beau texte [Be] de D. Bennequin.\par 
Je me propose d'exposer ici une premi\`ere approche, extr\^emement rudimentaire, de ce ''sens cognitif'' qui fait partie du programme de recherche de fondements cognitifs des math\'ematiques promu par le groupe ``G\'eom\'etrie et cognition'' form\'e par G. Longo, J. Petitot, et moi (voir [LPT]).\par
Il y a deux types d'ingr\'edients:\par\medskip\noindent
$\bullet$  Une interpr\'etation cognitive de certains objets math\'ematiques primitifs.\par\noindent Les progr\`es r\'ecents des neurosciences permettent de commencer \`a comprendre les bases biologiques du r\^ole constitutif irr\'eductible de l'espace et du temps dans notre repr\'esentation des ph\'enom\`enes, sur lequel ont insist\'e entre autres Kant, Poincar\'e, Hermann Weyl, Enriques.\par\medskip\noindent 
$\bullet$ Le moteur constitu\'e par ce que j'appelle la ''pens\'ee de bas niveau''.\par\noindent Ce vocable n'a ici rien de m\'eprisant, au contraire; il est inspir\'e par l'\'etude de la vision, qui n'a vraiment d\'ecoll\'e que lorsque les physiologistes ont essay\'e l'approche modeste consistant \`a tenter de mod\'eliser la ``vision de bas niveau'', qui est la partie ''imm\'ediate'' de la vision, avant la couleur et avant toute interpr\'etation (la neurophysiologie commence \`a donner une description de la mani\`ere dont cela se passe chez l'homme et du r\^ole des diff\'erentes aires visuelles du cortex).\par
 J'entends par l\`a les op\'erations de pens\'ee involontaires et tr\`es souvent {\it inconscientes}. Cela inclut des jugements involontaires comme celui de faire la distinction fixe/mobile, homog\`ene/inhomog\`ene, celui de ''comparer ce qui est comparable'', par exemple comparer la taille de deux objets de m\^eme nature, ou d\'etecter des r\'egularit\'es temporelles ou spatiales ou des sym\'etries, faire des analogies, faire la distinction entre un objet et ses attributs, ne pas distinguer des objets diff\'erents qui ont en commun des traits qui nous int\'eressent. Cela inclut aussi  des besoins fondamentaux (ou pulsions) de l'esprit humain comme la recherche obstin\'ee de causes ou d'origines, le fait de se demander si lorsque A implique B on a aussi B implique A, le besoin de se projeter dans l'avenir, de pr\'edire, de cr\'eer des rituels, ou encore celui de compl\'eter ce qui est incomplet, de d\'ecomposer un objet ou un m\'ecanisme complexe en objets ou m\'ecanismes simples.\par L'\'etude de cette pens\'ee de bas niveau reste \`a faire; elle est si proche de nous que nous ne la voyons pas. Elle appara\^\i t cependant \`a chaque \'etape de la pens\'ee, y compris dans les domaines les plus abstraits. Faire une liste raisonn\'ee des jugements et pulsions de bas niveau est en soi un d\'efi. Il me semble que ce sont Saint Augustin (en ce qui concerne le temps, dans {\it Les confessions}) puis de mani\`ere plus philosophique Thomas d'Aquin et Ma\"\i monide qui les premiers ont reconnu clairement, dans le cas particulier des attributs de la divinit\'e dans leurs religions respectives, que des automatismes de pens\'ee nous poussaient vers des paradoxes. ({\it cf}. [M]).\footnote{ Une des plus belles d\'ecouvertes de la pens\'ee religieuse, implicitement soulign\'ee dans le cadre du\break juda\"\i sme par Ma\"\i monide ({\it loc. cit.}), est celle de {\it l'\^etre qui ne fait rien d'autre qu'\^etre}, sans aucun attribut.} Dans un autre domaine, le d\'esir de bas niveau que toute proposition soit vraie ou fausse nous a aussi entra\^\i n\'es dans des orni\`eres en logique (voir le texte [G2] de Girard dans ce volume).\par\medskip
 L'id\'ee commune \`a ces deux ingr\'edients est que notre cerveau est le si\`ege d'activit\'es involontaires et inconscientes dont une partie ressemble \`a des math\'ematiques, et que ces activit\'es inconscientes interviennent dans nos activit\'es conscientes comme {\it r\'eservoir de sens}. \par Ce que nous percevons comme signification est en fait une r\'esonnance produite par notre physiologie entre notre pens\'ee consciente et la structure du monde telle que l'int\`egrent, de mani\`ere inconsciente, nos sens. Cette r\'esonnance est loin d'\^etre un isomorphisme au sens d\'ecrit au paragraphe suivant. Elle est au contraire assez souple pour pouvoir se propager le long de constructions formelles et langagi\`eres tr\`es \'elabor\'ees. \par\medskip\noindent
   Dans ce qui suit, je vais essayer d'illustrer cette id\'ee sur les exemples les plus simples.
\section{LA DROITE ET AUTRES LIEUX}\par\medskip\noindent
$\bullet$ 
{\bf La droite vestibulaire} (pour des renseignements pr\'ecis, voir le magnifique livre [B1] d'Alain Berthoz): le syst\`eme vestibulaire, situ\'e dans l'oreille interne, mesure avec une bonne pr\'ecision les acc\'el\'erations de toutes natures. Il joue un r\^ole important dans la survie des bip\`edes puisqu'il permet de r\'eagir tr\`es rapidement lorsque le sujet tr\'ebuche (la t\^ete part en avant avec une forte acc\'el\'eration). Il sert aussi \`a compenser automatiquement les mouvements de la t\^ete dus \`a la marche, ce qui nous permet de voir un monde stable pendant que nous marchons.  Mais surtout, avec la m\'emoire, c'est une centrale inertielle qui garde le souvenir de toutes nos acc\'el\'erations. D'apr\`es le principe de relativit\'e Galil\'eenne le seul mouvement qui ne donne aucun signal au syst\`eme vestibulaire (hormis le mouvement de la t\^ete, qui est bien s\'epar\'e du reste) est la marche \`a vitesse constante dans une direction fixe. La r\'ealit\'e neurophysiologique est bien plus riche que cette description sommaire.\par\noindent
Nous appellerons cet \'etat dynamique d'excitation minimale la {\it droite vestibulaire}.\par
Notons qu'elle est param\'etr\'ee par le temps, et scand\'ee par les pas. La seule mani\`ere de s'y rep\'erer est d'utiliser une horloge (par exemple notre horloge interne, ou celle donn\'ee par la marche).\par\medskip\noindent
$\bullet$ {\bf La droite visuelle} (voir [B1], cit\'e plus haut, et [N]): le nerf optique transmet les impulsions \'electriques produites dans la r\'etine par l'impact des photons \`a des neurones dont certains r\'eagissent \`a la pr\'esence dans une direction donn\'ee d'un petit segment ayant une orientation donn\'ee. Des neurones correspondant \`a des directions assez voisines s'excitent simultan\'ement s'ils d\'etectent la m\^eme orientation, bien plus que s'ils d\'etectent des orientations diff\'erentes. On pourrait dire que le ''transport parall\`ele'' est c\^abl\'e dans V1, mais la r\'ealit\'e est plus compliqu\'ee; la g\'eom\'etrie des connexions est loin de porter tout le fonctionnement, qui r\'esulte aussi d'interactions dynamiques entre les aires visuelles et pour cette raison on pr\'ef\`ere parler de {\it l'architecture fonctionnelle} des aires visuelles. Quoi qu'il en soit, cette architecture permet de d\'etecter des courbes et, parmi les courbes, d'isoler les DROITES, une droite \'etant ``une courbe qui a partout la m\^eme orientation'' selon la d\'efinition perceptuelle de J. Ninio dans [N]. Il semble que d\'ej\`a dans l'aire visuelle V1, par un r\'eseau d'excitations et d'inhibitions les lignes soient d\'etect\'ees. Les droites jouent un r\^ole sp\'ecial dans les aires visuelles situ\'ees en aval de V1\footnote{Je remercie Daniel Bennequin pour ses remarques sur cette partie.}. La d\'etection d'une droite (en fait, il s'agit d'un segment de droite) correspond donc \`a un \'etat dynamique d'excitation bien particulier (probablement extr\'emal en un certain sens) d'une assembl\'ee de neurones du cortex visuel. Nous appellerons un tel \'etat une {\it droite visuelle}. Il faut insister sur le fait que la possibilit\'e de l'existence de cette droite visuelle en neurophysiologie n'est pas du tout \'evidente et est un acquis relativement r\'ecent.\par\medskip\noindent
Il faut maintenant faire entrer en sc\`ene les relations tr\`es fortes existant dans le cerveau entre le syst\`eme visuel, le syst\`eme moteur et le syst\`eme vestibulaire. Les relations entre les trois ont \'et\'e \'etudi\'ees en particulier dans le laboratoire LPPA d'Alain Berthoz et sont si \'etroites ({\it cf}. [B1]) qu'elles justifient amplement l'intuition de Poincar\'e:\par\noindent Si les travaux cit\'es dans [B1] confirment l'intuition de Poincar\'e selon laquelle la position d'un objet dans l'espace visuel est reli\'ee par l'\'equivalent d'un changement de coordonn\'ees \`a l'ensemble des tensions musculaires correspondant au geste qu'il faut faire pour le saisir, il n'est pas abusif d'affirmer que l'\'evolution de nos syst\`emes de perception a aussi cr\'e\'e un isomorphisme entre la droite vestibulaire et la droite visuelle, dont j'ai soulign\'e dans [T3] l'importance pour les fondements de la signification de la droite en lui donnant le nom d'{\it isomorphisme de Poincar\'e-Berthoz}. \par
Ce n'est pas un isomorphisme au sens de la th\'eorie des ensembles munis de sructures: je ne peux pas exhiber une bijection d'un objet sur l'autre respectant la structure, et pour cause, puisque les deux droites sont des \'etats dynamiques d'assembl\'ees de neurones. Il vaudrait peut-\^etre mieux parler de ``correspondance'', mais la chose importante est que cette correspondance permet de transporter la structure de l'une des droites sur l'autre et une fois ce transport fait cette correspondance m\'erite qu'on lui applique le vocable ''isomorphisme''\footnote{Il s'agit plut\^ot d'un isomorphisme dans une cat\'egorie de percepts. Par ailleurs les connexions neuronales entre les aires motrices, visuelles et le syst\`eme vestibulaire sont -entre autres- le support mat\'eriel de l'isomorphisme mais sa nature est bien plus complexe.}.\par
Cet isomorphisme a en effet pour cons\'equence, par transport de structure de la droite vestibulaire sur la droite visuelle, de param\'etrer cette derni\`ere. C'est lui qui nous permet d'imaginer que nous avan\c cons le long d'une droite g\'eom\'etrique, c'est aussi lui qui nous permet d'accepter comme \'evidence que le temps est param\'etr\'e par une droite r\'eelle, ainsi que les propri\'et\'es de continuit\'e de l'ensemble des nombres r\'eels, ses propri\'et\'es archim\'ediennes (en un nombre fini de pas, je d\'epasse n'importe quel point), etc. C'est le genre d'images dont Einstein disait qu'elles jouaient un grand r\^ole dans sa pens\'ee. \par
Je sugg\`ere que notre droite math\'ematique a pour signification l'objet proto-\break math\'ematique obtenu par identification de la droite visuelle et de la droite vestibulaire, et que la partie de notre ''intuition math\'ematique'' qui est \`a l'oeuvre lors de d\'emonstrations portant sur les nombres r\'eels provient de cet objet proto-math\'ematique. Elle est l'origine de la construction des r\'eels comme compl\'etion des nombres rationnels ou de leur d\'efinition comme corps totalement ordonn\'e dans lequel tout sous-ensemble born\'e admet une borne sup\'erieure.\par Lorsque l'on me parle d'une suite $x_i$ de nombres r\'eels, je les ''vois'' sur la droite, et si l'on me dit qu'ils tendent vers $x$, je les vois osciller  peut-\^etre un peu autour de $x$ tout en s'en rapprochant in\'eluctablement. Bien s\^ur on peut exprimer tout cela avec des suites de rationnels ou des \'el\'ements d'un corps totalement ordonn\'e comme ci-dessus, mais alors l'intuition a bien du mal \`a suivre! Ceci est un exemple de r\'eponse au premier probl\`eme \'enonc\'e dans l'introduction\par
Un aspect int\'eressant de l'isomorphisme est que la droite visuelle et la droite vestibulaire n'ont pas les m\^emes automorphismes ''naturels''. Les automorphismes naturels de la droite vestibulaire sont en effet les seules translations $x\mapsto x+b$, ou tout au au plus les transformations affines $x\mapsto ax+b$ si l'on permet des mouvements uniformes de vitesses diff\'erentes. En ce qui concerne la droite visuelle, deux segments ferm\'es quelconques sont \'equivalents par un automorphisme affine de l'espace ambiant,  (qui appara\^\i t comme ``externe'' du point de vue vestibulaire ) et sont en particulier comparables au prix d'un automorphisme affine de l'espace ambiant \textit{respectant les longueurs} m\^eme si ils ne sont pas sur une m\^eme droite. Il est possible que ce simple fait, joint au besoin ''de bas niveau'' de comparer ce qui est comparable, ait conduit \`a la th\'eorie grecque du rapport ou logos, si magnifiquement expliqu\'ee dans le livre [F] de David Fowler. La question de comparer la longueur de deux segments est r\'esolue par les Grecs en les ramenant sur la m\^eme droite (par une isom\'etrie affine comme ci-dessus!) puis en regardant combien de fois (disons $a_0$) le petit ''va''dans le grand, puis combien de fois le segment qui reste du grand quand on a retir\'e $a_0$ fois le petit du grand ''va'' dans le petit, et ainsi de suite (antiphairesis). Le rapport de deux segments est ainsi d\'ecrit par une suit de nombres entiers, qui est infinie lorsque ce rapport est irrationnel. En termes plus modernes, c'est le d\'eveloppement en fractions continues du nombre r\'eel qui est le quotient des longueurs des deux segments mesur\'ees \`a l'aide d'une unit\'e de longueur quelconque. Une version d\'eg\'en\'er\'ee de cette op\'eration, sa restriction aux rapports qui sont des nombres rationnels, perdure dans notre enseignement sous le nom d'algorithme d'Euclide.   \par
Un autre aspect est que le continu visuel et le continu du mouvement se trouvent identifi\'es, ce qui a des cons\'equences \'enormes: cela permet d'inventer la notion de trajectoire param\'etr\'ee par le temps, donc finalement le concept de fonction, celui de vitesse (voir plus bas) et enfin l'espace-temps!. \par\medskip
Enfin, des d\'eveloppements math\'ematiques importants (et en particulier le cort\`ege d'id\'ees accompagnant l'ensemble de Cantor) doivent me semble-t-il leur naissance \`a la contradiction apparente entre la d\'efinition cognitive de la droite et le besoin (de bas niveau) de la comprendre comme ensemble de points. J'y reviendrai plus bas \`a propos des fronti\`eres. A ce propos il est fascinant de voir fonctionner la pens\'ee ''de bas niveau'' lorsque Aristote \'etudie l'hypoth\`ese que l'espace et le temps soient compos\'es d'indivisibles. \par
Le type d'analyse que je viens de d\'ecrire s'\'etend ensuite au plan et \`a l'espace. Par exemple le plan est param\'etrisable de deux mani\`eres du point de vue vestibulaire: coordonn\'ees cart\'esiennes ou coordonn\'ees polaires. Le plan visuel est celui de notre vision sans le relief, et il est possible que tout autre plan soit per\c cu comme un d\'eplacement de celui-ci. La description de l'espace visuel est compliqu\'ee: elle fait intervenir en particulier la vision binoculaire, et surtout le syst\`eme moteur; je renvoie \`a [B1].
 \par\medskip\noindent
$\bullet$ {\bf Ordinaux}\par\medskip\noindent 
La droite vestibulaire peut servir de support \`a une signification de l'infini~: ce qui n'est pas limit\'e, la marche sans fin sur la droite. Il faut lire les textes de Fontenelle cit\'es par Michel Blay dans [LAB], et admirer la finesse avec laquelle il d\'ecrit le malaise que suscite l'id\'ee de rajouter un nombre $\infty$ au bout de la suite des entiers. Ces textes r\'esument beaucoup de pages sur la diff\'erence entre l'infini ''potentiel'' et l'infini ''actuel'', dans lesquelles la r\'eflexion th\'eologique et la r\'eflexion math\'ematique se trouvent parfois extr\^emement proches. On voit, \`a partir de cette \'epoque, les math\'ematiciens accepter progressivement l'infini dans leur bo\^\i te \`a outils et en pr\'eciser la signification.\par\noindent Par exemple d\`es la premi\`eme moiti\'e du 19$^{\hbox{\rm \`eme}}$ si\`ecle Bolzano avait r\'ealis\'e l'importance du fait qu'un ensemble infini peut \^etre en bijection avec un de ses sous-ensembles stricts (par exemple les entiers et les entiers pairs). La perc\'ee devra cependant attendre plus d'un demi-si\`ecle l'\'epoque de Cantor.\par  Par contre, notre intuition du monde accepte beaucoup plus difficilement l'id\'ee de marcher {\it depuis} un temps infini.  C'est peut-\^etre de l\`a que provient notre pr\'ef\'erence spontan\'ee pour les bons ordres, ceux dans lesquels il ne peut exister de suite infinie d\'ecroissante (vers un point d'o\`u l'on ne peut pas revenir en un nombre fini de pas!), en m\^eme temps que notre besoin (de bas niveau) si fort que tout ait une origine ou une cause (revoici la th\'eologie!).\par L'axiome du choix qui \'enonce que {\it l'on peut choisir simultan\'ement un \'el\'ement dans chacun des ensembles d'une collection \'eventuellement infinie d'ensembles non vides} et le fait que {\it tout ensemble peut \^etre bien ordonn\'e}, c'est \`a dire muni d'un ordre {\it total} (ce qui signifie que deux \'el\'ements quelconques sont comparables) sans suite d\'ecroissante  infinie, correspondent de la mani\`ere la plus \'evidente \`a deux pulsions de bas niveau: le besoin que tout ait une origine et celui, peut-\^etre de moins bas niveau et plus math\'ematique, de pouvoir choisir ind\'efiniment. Ceci est un exemple de r\'eponse \`a la seconde question pos\'ee dans l'introduction.\par Ces deux assertions sont {\it a priori} diff\'erentes mais l'\'equivalence des deux \'enonc\'es math\'ematiques est un th\'eor\`eme de th\'eorie des ensembles. On peut consid\'erer que cela  rend manifeste une interaction entre les deux ''\'evidences'' du cas fini que ces \'enonc\'es \'etendent au cas infini.\par
Un grand moment de la pens\'ee math\'ematique est la cr\'eation par Cantor de la th\'eorie des ensembles bien ordonn\'es ou du moins totalement ordonn\'es, et en particulier la prise de conscience du fait qu'il y en avait beaucoup. Les classes d'isomorphisme (= bijection croissante) d'ensembles bien ordonn\'es sont appel\'ees {\it ordinaux}. Les ordinaux peuvent eux-m\^emes \^etre munis d'un ordre: $\alpha \leq \alpha'$ signifie qu'il existe une injection croissante d'un repr\'esentant de $\alpha$ dans un repr\'esentant de $\alpha '$. Cet ordre fait des ordinaux un ensemble\footnote{Note technique: pour la plupart des math\'ematiciens, qui travaillent dans un mod\`ele fixe de la th\'eorie des ensembles, con\c cu comme un univers (mental), les ordinaux constituent un ensemble, alors qu'en logique leur statut est plus complexe.}  bien ordonn\'e. Chacun des ordinaux finis est repr\'esent\'e par un ensemble d'entiers $\{1,\ldots ,n\}$ muni de l'ordre usuel, et on peut le d\'esigner par l'entier $n$. Le premier exemple infini est la classe $\omega$ de l'ensemble ordonn\'e $\N$ des entiers.  C'est par la consid\'eration des ordinaux sup\'erieurs \`a $\omega$ que l'on peut se repr\'esenter $\infty$, ou plut\^ot $\omega$, d'une mani\`ere que notre intuition accepte facilement: c'est {\it le plus petit} ordinal plus grand que tous les ordinaux qui sont les classes d'ensembles finis!.  \par\medskip\noindent
 Plus pr\'ecis\'ement, consid\'erons l'ensemble $\Z^2$ des paires d'entiers positifs ou n\'egatifs, ordonn\'e lexicographiquement. L'in\'egalit\'e $(i,j)< (i',j')$ signifie que $i$ est inf\'erieur \`a $i'$ ou que $i=i'$ et $j$ est inf\'erieur \`a $j'$. C'est bien l'\'equivalent num\'erique de l'ordre du dictionnaire. \par C'est bien un ordre total, mais il ne fait pas de $\Z^2$ un ensemble bien ordonn\'e. Un sous-ensemble bien ordonn\'e de $\Z^2$ consiste en points $(i,j)$ dont la premi\`ere coordonn\'ee est born\'ee inf\'erieurement (disons $i\geq 0$ pour fixer les id\'ees) et tel que de plus, {\it pour chaque $i$}, la seconde coordonn\'ee  des points de l'ensemble dont la premi\`ere coordonn\'ee vaut $i$ est born\'ee inf\'erieurement. Par exemple l'ensemble $\N^2$ des points dont les deux coordonn\'ees sont positives est bien ordonn\'e. L'ensemble des points dont la premi\`ere coordonn\'ee est positive ne l'est pas. Consid\'erons le sous-ensemble bien ordonn\'e de $\Z^2$ obtenu comme ceci: pour chaque entier $i\geq 0$ choisissons un entier positif ou n\'egatif $\beta(i)$. Notre ensemble est:
  $$E_\beta=\{(i,j)\in \Z^2\vert i\geq 0,\ \ j\geq \beta(i)\}.$$ Nous pouvons maintenant faire en pens\'ee le parcours suivant, qui \'etend la droite vestibulaire:\par\noindent
  Partons du point $(0,\beta (0))$ et augmentons la seconde coordonn\'ee d'une unit\'e \`a la fois. Nous marchons sur la droite vestibulaire identifi\'ee \`a l'axe ``vertical'' des points dont la premi\'ere coordonn\'ee vaut $0$. Apr\`es une infinit\'e de pas le point de l'ensemble $E_\beta$ qui se pr\'esente \`a nous naturellement est le point $(1,\beta (1))$, qui est un parfait repr\'esentant de l'infinit\'e de pas que nous venons de faire; c'est le plus petit \'el\'ement de notre ensemble $E_\beta$ qui soit plus grand que tous les entiers, identifi\'es aux points de coordonn\'ee $(0,j); j\geq \beta(0)$. Appelons le $\omega$. Puis poursuivons notre marche, cette fois-ci en montant parmi les points de $E_\beta$ dont la premi\`ere coordonn\'ee vaut $1$. Le point $(1,\beta (1)+1)$ s'appellera $\omega +1$, et ainsi de suite. Le point $(2,\beta (2))$ s'appellera naturellement\footnote{En fait, il s'appellerait tout aussi naturellement $2\omega$ mais la multiplication des ordinaux au moyen de produits lexicographiques n'est pas commutative et avec la convention usuelle $2\omega$ est \'egal \`a $\omega$.}$\omega 2$, et lorsque nous aurons \'epuis\'e toutes les valeurs de $i$ nout atteindrons $\omega^2$, qui est l'ordinal repr\'esent\'e par $\N^2$ muni de l'ordre lexicographique, mais pour repr\'esenter cet ordinal par un point, nous devrons nous placer dans $\Z^3$ muni de l'ordre lexicographique. Nous pourrions \'evidemment ajouter formellement \`a l'ensemble $\N$ (ou \`a $\N^2$) un \'el\'ement $\omega$ (ou $\omega^2$) et \'etendre la relation d'ordre en d\'ecidant que ce nouvel \'el\'ement est plus grand que tous les \'el\'ements de $\N$ (ou $\N^2$). Cette construction brutale n'a rien de g\'eom\'etrique et reste proche de ce que Ch\^atelet appelle je crois ''la m\'ediocrit\'e de l'agr\'egat''. La repr\'esentation par les produits lexicographiques est bien plus riche et m\`ene en particulier \`a la d\'efinition du {\it produit} des ordinaux et \`a la description de la structure de "l'ensemble" des ordinaux. \par\medskip\noindent Nous avons ainsi satisfait une pulsion de la pens\'ee de bas niveau: compl\'eter ce qui est incomplet, en plongeant les entiers dans un ensemble bien ordonn\'e plus grand. Celui-ci poss\`ede \`a son tour un infini, que l'on peut r\'ealiser dans un ensemble bien ordonn\'e plus grand, et ainsi de suite\footnote{Rappelons que cette soudaine multiplication des infinis n'a pas \'et\'e accept\'ee sans difficult\'es par des math\'ematiciens proches de la physique dont  l'intuition rechignait. Mais tous ont fini par se rallier au point de vue de Hilbert affirmant que ''nul  ne le chasserait du paradis cr\'e\'e par Cantor''. Aujourd'hui encore, les anglo-saxons utilisent la terminologie rassurante de {\it linearly ordered set} pour {\it ensemble totalement ordonn\'e}.}. L'ordinal $\omega$ appara\^\i t d\'esormais comme un nombre assez semblable aux autres; il en diff\`ere par le fait qu'il faut une infinit\'e de pas pour l'atteindre, et que l'on ne peut pas parler de $\omega -1$ puisqu'il n'a pas de pr\'ed\'ecesseur.\par\noindent Pour rassurer le lecteur, ou plut\^ot pour satisfaire sa pulsion de bas niveau concernant l'unicit\'e des repr\'esentations, j'ajouterai que cette repr\'esentation est, \`a un isomorphisme d'ensembles ordonn\'es (= bijection croissante) pr\`es, ind\'ependante du choix de la fonction $i\mapsto \beta(i)$ de $\N$ dans $\Z$. L'ensemble bien ordonn\'e $\N^2$ correspond au choix $\beta(i)=0$ pour tout entier $i$.\par\par\medskip\noindent
$\bullet$ {\bf Bords et fronti\`eres}\par\medskip\noindent
La droite visuelle et la droite vestibulaire portent toutes deux des fronti\`eres, qui sont le {\it bord} ou l'{\it extr\'emit\'e} pour l'une, et la {\it fin} (du mouvement) pour l'autre. C'est une propri\'et\'e fondamentale du continu, d\'ej\`a rep\'er\'ee par Aristote, que d'\^etre constitu\'e de fronti\`eres. La correspondance entre une fronti\`ere spatiale et une fronti\`ere temporelle dans l'isomorphisme de Poincar\'e-Berthoz a de nombreuses cons\'equences. Je veux me concentrer sur le fait que le mouvement pouvant s'interrompre \`a tout instant, la droite visuelle, par transport de structure, devient divisible comme le temps et que cette division est manifest\'ee par le choix d'un point d'arr\^et.\par  En fait la r\'eflexion classique sur la divisibilit\'e de l'espace et du temps a je pense deux versions qui sont en quelque sorte duales: le continu n'est pas ''form\'e'' de points, c'est \`a dire qu'il n'est pas seulement un ensemble de points, mais il contient des points, qui y apparaissent comme des {\it coupures} qui le s\'eparent en morceaux adjacents, et qui sont pr\'ecis\'ement les coupures\footnote{Une coupure est une partition de l'ensemble des nombres rationnels (ou r\'eels) en deux sous ensembles non vides tels que tout \'el\'ement du premier soit inf\'erieur \`a tout \'el\'ement du second. Par exemple les rationnels inf\'erieurs \`a $\sqrt 2$ et ceux qui sont sup\'erieurs. Une coupure des rationnels d\'etermine un unique nombre r\'eel.} de Dedekind. L'avatar temporel de la coupure est {\it l'instant}. L'observation fondamentale de Dedekind (voir [De]) est que, comme le temps, le continu unidimensionel contient comme points {\it toutes} les fronti\`eres (ou {\it coupures}) que l'on peut y d\'efinir et rien d'autre, ce qui n'est pas le cas de l'ensemble des nombres rationnels contenus dans un intervalle donn\'e. D'autre part un continuum (compact) de nombres r\'eels, c'est \`a dire  un intervalle $[a,b]$ ferm\'e (= contenant ses extr\'emit\'es), poss\`ede une structure {\it autosimilaire}: il est r\'eunion de sous-intervalles stricts qui lui sont isomorphes, et m\^eme par des isomorphismes naturels. De m\^eme que l'on peut it\'erer ind\'efiniment la s\'eparation en deux par le choix d'un point, on peut it\'erer \`a l'infini la d\'ecomposition correspondante d'un intervalle ferm\'e comme r\'eunion de deux intervalles ferm\'es plus petits. Cela exprime l'homog\'eneit\'e du continu. La construction de l'ensemble de Cantor et celle de ses analogues en dimension sup\'erieure reposent sur cette autosimilarit\'e de la droite, du plan, de l'espace. En fait, la construction de Cantor est la construction la plus simple de ce point de vue pour un sous-ensemble de l'intervalle $[0,1]$ qui soit totalement discontinu, c'est \`a dire d\'epourvu de fronti\`eres (techniquement: ne contenant aucun intervalle), et non d\'enombrable ainsi que son compl\'ementaire. Tant la notion de fronti\`ere que celle d'autosimilarit\'e sont essentiellement topologiques. La diff\'erence entre points irrationnels et points rationnels dispara\^\i t lorsqu'on les consid\`ere en tant que coupures, et il ne reste que des bords tous isomorphes, ce qui exprime aussi  l'homog\'eneit\'e du continu. \par Les paradoxes de Z\'enon sur l'impossibilit\'e du mouvement reposent d'une part sur la question de la finitude d'une somme infinie de termes, qui heurte notre intuition vestibulaire li\'ee \`a la marche, et d'autre part sur la question d'atteindre ou non en un temps fini une fronti\`ere situ\'ee \`a distance finie. En l'absence du concept de vitesse la question est difficile! Or, \`a partir de l'antiphairesis, il a fallu du temps \`a la pens\'ee scientifique pour accepter de diviser des longueurs par des temps (rappelez-vous: comparer ce qui est comparable). On peut dire que la d\'efinition de la vitesse est une retomb\'ee num\'erique de l'isomorphisme de Poincar\'e-Berthoz, qui n'\'etait possible qu'apr\`es que l'antiphairesis ait c\'ed\'e la place au quotient des longueurs pour comparer deux segments.\footnote{Les diagrammes d'Oresme du chapitre 2 de [C] peuvent s'interpr\'eter comme une tentative de pr\'esenter ''\`a la grecque'' (c'est \`a dire que les produits sont vus comme des aires) au moyen de diagrammes la longueur parcourue lors d'un mouvement comme un {\it produit} de ''quelque chose'' par le temps, ce qui \'etait sans doute la seule mani\`ere possible \`a l'\'epoque de diviser de mani\`ere structurelle (par opposition \`a num\'erique) une longueur par un temps pour \textit{d\'efinir} la vitesse moyenne.} Cela souligne le fait que notre espace mental ne nous incite pas \`a nous demander combien de fois un temps ``va'' dans une longueur, m\^eme si l'on peut mesurer des distances en heures ou journ\'ees de voyage. Il nous faut d'abord les transformer en nombres par une mesure. La notion qualitative de vitesse fait \'evidemment partie de notre exp\'erience primitive du monde et l'animal poursuivant sa proie en a une perception extr\^emement pr\'ecise, mais sa d\'efinition quantitative est assez complexe.\par\noindent
\section{L'APPROCHE COGNITIVE N'EST PAS ANTHROPOMORPHIQUE}\par\noindent

\par Rappelons les consid\'erations de Poincar\'e dans [P] sur la ``M\'ecanique anthropomorphique'', \`a propos des notions de force et de chaleur dont nous avons ''l'intuition directe''. Je renvoie aussi \`a ce qui vient d'\^etre dit sur la vitesse:\par\noindent
{\it ...cette notion subjective ne peut se traduire en nombre, donc elle ne sert \`a rien...}\par\noindent
et, plus loin:
\par\noindent
{\it L'anthropomorphisme a jou\'e un r\^ole historique consid\'erable dans la gen\`ese de la m\'ecanique; peut-\^etre fournira-t-il encore quelques fois un symbole qui para\^\i tra commode \`a quelques esprits; mais il ne peut rien fonder qui ait un caract\`ere vraiment scientifique, ou un caract\`ere vraiment philosophique}.\par\noindent 
Mais il y a une grande diff\'erence entre le fait de rechercher les fondements de la science dans les sensations et le fait de rechercher les fondements de la signification de la science dans la structure des syst\`emes de perception et d'action et dans la pens\'ee de bas niveau. Ce que nous proposons ici diff\`ere beaucoup d'une approche anthropomorphique au sens o\`u l'entend Poincar\'e. \par Celui-ci r\'efute \`a juste titre l'id\'ee d'une relation directe entre nos sensations et la science quantitative. Il d\'enie m\^eme, encore \`a juste titre, une signification philosophique \`a ces relations mal d\'efinies entre la description scientifique des ph\'enom\`enes et la perception qualitative que nous pouvons en avoir.\par Je propose l'id\'ee d'une relation tr\`es forte entre les {\it modes d'int\'egration} des donn\'ees de nos diff\'erents syst\`emes de perception et de leurs \'etats qui permettent une perception unifi\'ee de notre environnement et la {\it signification} des objets math\'ematiques que nous utilisons pour d\'ecrire cet environnement.
\par J'y ajoute l'importance pour la construction scientifique des pulsions simplificatrices inconscientes que j'appelle de bas niveau, auxquelles Poincar\'e ne fait pas r\'ef\'erence mais dont la reconnaissance me semble aussi n\'ecessaire que celle du caract\`ere quantitatif sur lequel il insiste.\par
En revanche, l'approche cognitiviste met en \'evidence le fait que nos constructions math\'ematiques sont fortement d\'ependantes d'un syst\`eme perceptif et de pulsions que nous partageons en grande partie avec les primates ({\it cf.} [T2], [T3]). Des \^etres pensants ayant les perceptions de poulpes construiraient peut-\^etre les m\^emes objets math\'ematiques mais pas pour les m\^emes raisons ni avec la m\^eme signification. 
\section{L'APPROCHE COGNITIVE N'EST PAS REDUCTIONNISTE}\par\noindent

\par Il faut peut-\^etre d\'efendre ce point de vue contre une accusation de r\'eductionnisme. Bien s\^ur la signification n'est pas ``c\^abl\'ee'' dans notre cerveau et n'est pas r\'eductible \`a la physico-chimie. L'id\'ee qui me semble la plus pertinente est encore due \`a Thom: L'investissement des ''pr\'egnances'' de sens sur des ''saillances'' de langage (ici la ''droite math\'ematique'') est une op\'eration du type ''cr\'eation de langage'' qui n'ob\'eit pas au principe de Curie (les sym\'etries des causes se retrouvent dans les sym\'etries des effets ) et par cons\'equent n'est pas susceptible d'une description formelle ({\it cf.} [Th], pp. 113-116 et 310-311). S'il y a r\'eduction, elle est indescriptible! \par Une autre mani\`ere de le dire est que les neurosciences nous apprennent entre autres qu'au niveau physico-chimique notre cerveau est le si\`ege de ph\'enom\`enes d'une complexit\'e dont notre ``pens\'ee pensante'' n'a aucune id\'ee. Il y a des dynamiques sur plusieurs \'echelles qui interagissent, et bien d'autres sources de complexit\'e.  Cela \'etant admis, quelle signification aurait une ''r\'eduction'' de la pens\'ee \`a des ph\'enom\`enes physico-chimiques? Ce n'est vraiment pas du tout comparable avec la r\'eduction de la chimie aux interactions mol\'eculaires, qui d'ailleurs est loin d'\^etre compl\`etement faisable en pratique parce qu'elle est elle-m\^eme d'une grande complexit\'e. La biologie n'a pas encore trouv\'e son G\"odel, celui qui convaincra les pratiquants de la discipline que tout ce qui est observable n'est pas rationnellement d\'eductible de principes et de faits \'el\'ementaires.  En revanche, l'analyse neurophysiologique,  l'imagerie et la mod\'elisation apportent des informations pr\'ecieuses (et parfois d'une pr\'ecision spectaculaire, voir [Pe2]) sur des m\'ecanismes pr\'ecis \`a des \'echelles pr\'ecises, et sugg\`erent de nouveaux concepts.\section{IT'S THE GEOMETRY, STUPID!}\par\noindent{\small (Voir [Car])}
\par\medskip
La seule signification de la phrase kantienne ``L'espace est une donn\'ee {\it a priori} de la conscience'' qui me soit accessible est que notre rapport \`a l'espace est {\it imm\'ediat} au sens \'etymologique.\par\noindent C'est me semble-t-il d\^u exactement \`a cette identification de nos perceptions visuelles, vestibulaires et motrices d\'ecrite dans [B1], et c'est une source in\'epuisable d'images et d'analogies, qui nous sont des ''mod\`eles'' proto-math\'ematiques de dynamiques tr\`es \'elabor\'ees. \par Il faudrait \'etudier la mani\`ere dont l'isomorphisme de Poincar\'e-Berthoz s'\'etend aux courbes, trajectoires dynamiques aussi bien que visuelles, qui ont les m\^emes propri\'et\'es topologiques que la droite, mais qui ne sont pas aussi \'evidemment autosimilaires. \par
Le concept de trajectoire ram\`ene (au prix de la construction d'un espace de phases) notre conception de l'\'evolution temporelle de tout syst\`eme m\'ecanique, si compliqu\'e soit-il, \`a la trajectoire d'un caillou que l'on lance dans {\it notre} espace. D'ailleurs il serait sans doute plus pertinent, comme l'avait d\'ej\`a vu Poincar\'e, de dire que l'espace-temps, ou plut\^ot le {\it mouvement} est une donn\'ee {\it a priori} de la conscience.\par\noindent L'int\'egration du temps comme dimension est permise par l'isomorphisme de Poincar\'e-Berthoz, mais celui-ci n'\'epuise pas le sujet; d'autres ph\'enom\`enes perceptifs de nature analogue r\'egissent sans doute la cr\'eation de l'id\'ealit\'e math\'ematique du mouvement. \par Une partie des math\'ematiques consiste \`a cr\'eer, au moyen de la m\'ethode axiomatique, des espaces dans lesquels nous pouvons nous mouvoir -en pens\'ee- ''presque'' comme dans notre espace usuel, mais dont les points ont des significations tr\`es abstraites: ils peuvent repr\'esenter des fonctions, des ensembles de ''nombres'', des mesures, des lois de probabilit\'e, des objets g\'eom\'etriques, etc. Nous y transportons la notion de distance \`a laquelle nous sommes habitu\'es depuis des centaines de g\'en\'erations, ou encore celle des sym\'etries de l'espace, et elle permet de d\'emontrer par exemple des r\'esultats portant sur les solutions d'\'equations aux d\'eriv\'ees partielles ou, dans le cas de la G\'eom\'etrie des nombres de Minkowski, sur les solutions d'\'equations diophantiennes. \par J'aime r\'esumer cela en disant que nous ne comprenons un th\'eor\`eme que lorsque nous avons r\'eussi \`a l'expliquer au primate qui est en nous (et qui a \'et\'e entra\^\i n\'e \`a cela pendant notre formation).\par Mais presque tout reste \`a comprendre de la mani\`ere dont la repr\'esentation formelle int\`egre ce qui provient de l'utilisation dans ces espaces de notre intuition de l'espace \`a trois dimensions ou plut\^ot, encore une fois, du mouvement. Cela revient \`a d\'ecrire une toute petite partie de la mani\`ere dont l'inconscient donne au langage sa structure et sa signification ({\it cf.} [T2]). Je renvoie aussi au texte [G1] de Girard. \par\smallskip\noindent
Par ailleurs l'interpr\'etation perceptive de certains objets fondamentaux comme l'ensemble des entiers n'est certainement pas \'epuis\'ee par une approche g\'eom\'etrique, et la marche n'est qu'un exemples des actions ou ph\'enom\`enes r\'ep\'etitifs ou p\'eriodiques dont le rapprochement avec le d\'enombrement des objets fonde me semble-t-il la signification des entiers. La danse trouve probablement en partie son origine dans un lien perceptif entre la marche et d'autres mouvements p\'eriodiques du corps et les perceptions auditives p\'eriodiques.\par
Il semble par ailleurs que, au cours du d\'eveloppement, notre perception imm\'ediate des nombres (ou plut\^ot des cardinaux finis) est beaucoup plus affect\'ee par l'acquisition du langage que notre perception de l'espace. 

\section{CONCLUSION:}\par\medskip\noindent

Il me semble que l'on peut dire que l'isomorphisme de Poincar\'e-Berthoz {\it constitue} l'objectivit\'e de la signification de l'id\'ealit\'e math\'ematique appel\'ee droite. Cette signification n'a donc pas d'existence absolue, ind\'ependante d'un \^etre dont la perception puisse construire cet isomorphisme, et son objectivit\'e est plut\^ot un fait de nature intersubjective entre de tels \^etres.\par 
La construction classique de la droite r\'eelle \`a partir des entiers de la th\'eorie des ensembles, par passage aux nombres rationnels puis par compl\'etion (ou construction de coupures) jusqu'aux nombres r\'eels {\it constitue} l'objectivit\'e de la v\'erit\'e des \'enonc\'es concernant l'id\'ealit\'e math\'ematique appel\'ee droite. Si l'on admet l'objectivit\'e de l'ensemble des entiers, celle de la droite est inattaquable. Cependant, des \^etres ayant d'autres perceptions la constitueraient peut-\^etre par un autre chemin. \par
La relation entre ces deux constructions d'id\'ealit\'es  n'est pas imm\'ediatement \'evidente. L'usage interm\'ediaire des rationnels peut \^etre vu comme une premi\`ere tentative d'atteindre la divisibilit\'e du continu dont il a \'et\'e question plus haut. Son \'echec est manifest\'e par l'existence de coupures irrationnelles dans l'ensemble des nombres rationnels et l'ajout de celles-ci suffit pour obtenir l'objet qui, comme la droite visuelle, contient toutes les fronti\`eres que l'on peut imaginer d'y tracer. Le fait d'admettre que ces deux constructions ''parlent de la m\^eme chose'' est entre autres une assertion sur le rapport entre discret et continu et est un vrai sujet de philosophie des math\'ematiques, dont l'interpr\'etation donn\'ee ici me semble compatible avec (mais pas du tout \'epuis\'ee par) le Platonisme G\"odelien dont Petitot se fait l'avocat dans [LAB].\par
Les r\'eflexions pr\'esent\'ees ici partagent avec celles de Gilles Ch\^atelet dans son superbe livre [C] une m\^eme qu\^ete ''irrationnelle'' du sens des constructions physiques ou math\'ematiques mais les chemins suivis, s'ils se rapprochent parfois, ne sont pas les m\^emes. Bien que les r\'eflexions de Ch\^atelet soient beaucoup plus \'elabor\'ees que celles de ce texte, j'esp\`ere pouvoir faire ailleurs quelques comparaisons.\par
J'esp\`ere aussi avoir convaincu le lecteur par ces quelques exemples encore tr\`es peu d\'evelopp\'es qu'il \'etait int\'eressant de rechercher dans la structure extr\^emement riche de notre syst\`eme inconscient de perception de notre environnement et de pulsions simplificatrices le fondement de la signification de certaines des id\'ealit\'es et de certains des r\'esultats des math\'ematiques.\par\smallskip
{\small \it Je remercie Philippe Courr\`ege de sa lecture aig\"ue et de ses critiques constructives.}
\par\vskip.5truecm\noindent
\centerline{Bibliographie}\par\medskip\noindent
[B1] Berthoz, A., 1997, {\it Le sens du Mouvement}, Editions Odile Jacob.\par\noindent
[B2] Berthoz, A., 2003, {\it La D\'ecision}, Editions Odile Jacob.\par\noindent
[Be] Bennequin, D., 1994, {\it Questions de Physique galoisienne}, in ''Passion des formes'', Mich\`ele Porte, Editeur, Presses de l'ENS Fontenay, diffusion Ophrys. \par\noindent
[B-L] Bailly, F., Longo, G.,  2006, {\it Math\'ematiques et sciences de la nature, la singularit\'e physique du vivant}, Hermann, Paris.\par\noindent
[C] Ch\^atelet, G., 1993, {\it Les enjeux du mobile, Math\'ematique, Physique, Philosophie}, Coll. Des Travaux, Le Seuil.\par\noindent
[Ch] Changeux, J.-P., 1983, {\it L'homme neuronal}, Ed. Odile Jacob.\par\noindent
[C-C] Connes, A. et Changeux, J.-P., 2000, {\it Mati\`ere \`a penser}, Ed. Odile Jacob.\par\noindent
[Car] Carville, J., 1992, auteur du memento ''It's the economy, stupid!'' de la campagne pr\'esidentielle de Bill Clinton.\par\noindent
[D] Damasio, A., 1995, {\it L'erreur de Descartes}, Ed. Odile Jacob.\par\noindent
[De] Dedekind, R., 1872, Trait\'es sur la Th\'eorie des Nombres, trad. C. Duverney, Pr\'eface de G. Wanner, Editions du Tricorne, Gen\`eve 2006.\par\noindent
[F] Fowler, D, 1987, {\it The Mathematics of Plato's Academy}, Clarendon Press, Oxford.\par\noindent
[G1] Girard, J.-Y., 2003, {\it La logique comme g\'eom\'etrie du cognitif}, pr\'epublication, voir: http://iml.univ-mrs.fr/~girard/Articles.html.\par\noindent
[G2] Girard, J.-Y., 2007 {\it De la syllogistique \`a l'iconoclasme}, ce volume.\par\noindent
[G3] Girard, J.-Y., 2000 {\it Les Fondements des Math\'ematiques}, 169$^{\hbox{\rm \`eme}}$ conf\'erence \`a l'Universit\'e de tous les savoirs.\par\noindent
[LAB] {\it Le labyrinthe du continu}, Colloque de Cerisy, J-M Salanskis, H. Sinaceur, Editeurs, Springer Verlag Paris 1992.\par\noindent
[LPT] Longo, G., Petitot, J., Teissier, B., 1999, voir ``Motivations g\'en\'erales'',  in ``G\'eom\'etrie et cognition'', on http://www.di.ens.fr/users/longo/geocogni.html\par\noindent
[M] Ma\"\i monide, {\it Le guide des \'egar\'es}, Editions Verdier, 1979.\par\noindent
[N] Ninio, J., 1989, {\it L'empreinte des sens}, Ed. du Seuil.\par\noindent
[P] Poincar\'e, H., {\it La science et l'hypoth\`ese}, Flammarion, Paris.\par\noindent
[Pe1] Petitot, J., 1985, {\it Morphog\'en\`ese du sens}, PUF, Paris.\par\noindent
[Pe2] Petitot, J., 2006, {\it Neurog\'eom\'etrie des architectures fonctionnelles de la vision}, Journ\'ee annuelle de la Soci\'et\'e Math\'ematique de France, Juin 2006. Publ. SMF, Paris. Disponible sur http://smf.emath.fr/VieSociete/JourneeAnnuelle/2006/\par\noindent
[T1]  Teissier, B., 1994, {\it Des mod\`eles de la Morphog\'en\`ese \`a la Morphog\'en\`ese des mod\`eles}, in ``Passion des formes'', Mich\`ele Porte, coordonnateur, ENS Editions Fontenay-Saint-Cloud, diffusion Ophrys.\par\noindent
[T2]  Teissier, B., 2004, {\it Le mur du langage}, in ``Le r\'eel en math\'ematiques, Math\'ematiques et psychanalyse'',  P. Cartier et N. Charraud, \'editeurs, Editions Agalma, diffusion Le Seuil.\par\noindent
[T3] Teissier, B., 2005, {\it Protomathematics, perception and the meaning of mathematical objects}, in "Images and Reasoning", edited by P. Grialou, G. Longo, M. Okada, CIRM, Keio University, Tokyo 2005.\par\noindent
[Th] Thom, R., 1972, Stabilit\'e structurelle et Morphog\'en\`ese, W.A. Benjamin, Inc., Reading, Massachusetts, Inter\'edition, Paris.

\end{document}